\documentclass[11pt,a4paper]{article}

\usepackage{titletoc,amsmath,mathrsfs,amsfonts,amssymb}
\usepackage{txfonts}
\usepackage{color}

\makeatletter
\def\@evenhead{
\vbox{\hbox to \textwidth {}{\hspace{0mm}{\footnotesize
\thepage}}{\hspace{10cm} {\footnotesize {Ling-di Wang}}} \protect\vspace{1truemm}\relax \hrule depth0pt
height0.15truemm width\textwidth}}
\def\@evenfoot{}
\def\@oddhead{\vbox{\hbox to \textwidth
{{\hspace{0cm}{\footnotesize Mixed eigenvalues of discrete $p\,$-Laplacian}\hfill{\footnotesize
\thepage}}\hspace{0mm}}{} \protect\vspace{1truemm}\relax\hrule
depth0pt height0.15truemm width\textwidth}}
\def\@oddfoot{}
\makeatother

\def\scr{\mathscr}
\newfont{\aaa}{cmb10 at 19pt}
\newfont{\bbb}{cmb10 at 11pt}
\newtheorem{thm}{\bf Theorem}[section]
\newtheorem{lem}[thm]{\bf Lemma}
\newtheorem{exam}[thm]{\bf Example}

\def\minus{\setminus}

\begin{document}

\thispagestyle{empty}
\noindent{\aaa{Mixed  eigenvalues of {\LARGE $\pmb p$\,}-Laplacian on trees}}\\[1mm]

\noindent{\bbb Ling-Di WANG}\\[-1mm]

\noindent\footnotesize{ School of Mathematics and Statistics, Henan University, Kaifeng 475001, China}\\[6mm]

\noindent{\bf Abstract} {The purpose of the paper is to present quantitative estimates for the principal eigenvalue of discrete p-Laplacian on the set of rooted trees. Alternatively, it is studying the optimal constant of a class of weighted Hardy inequality. Three kinds of variational formulas in different formulation for the mixed principal eigenvalue of  p-Laplacian  on the set of trees with unique root as Dirichlet boundary are presented. As their applications,  we obtain a basic estimate of  the  eigenvalue on trees.
}

\section{ Introduction}\label{c1}
In \cite{r1, W-Z}, mixed principal eigenvalue for birth-death process on line were studied. Inspired by analogies research for that, mixed principal p-Laplacian on line were studied  in \cite{CWZ,CWZ1}.  We shall extend the related results to a more general setting, investigating the quantitative estimates of the mixed principal p-Laplacian on trees with unique root as Dirichlet  boundary. A basic result on the property of eigenvalue of p-Laplacian on trees, which is a key point for the extension, will be  obtained.

By a tree, denote $T$, we mean $T$ a undirected, connected, locally finite graph without cycles. One distinguished vertex, say $o$, is called the root. For any vertex $i$, the number of edges on the unique simple path between $i$ and the root $o$ is called the level of $i$ and denote $|i|$. Let $E$ be the edge set and $V$ be the vertexes of $T$. The vertexes at level $|i|+1$ (correspondingly, $|i|-1$) that are adjacent to $i$ are called children (correspondingly,  parents) of $i$. Throughout the paper, we assume that trees are locally finite (i.e., each vertex has finite chilren).

To be specified, $J(i)$ is the set of children of vertex $i$ and $i^*$ is a parent of $i$. Operator $\Omega_p$ we focusing on in the paper is of the form
 $$\aligned\Omega_p f(i)&=\sum_{j\in J(i)}\nu_j|f_j-f_i|^{p-2} (f_j-f_i)+\nu_i|f_{i^*}-f_i|^{p-2}(f_{i^*}-f_i), \qquad i\in V\minus\{o\},\endaligned$$
where $\{\nu_i: i\in V\}$ is a positive sequence.
 We concentrate on estimating the $p$-Laplacian eigenvalue on a tree, which is described as follows:
 \begin{eqnarray}
\text{``Eigenequation'':}& \Omega_pg(k)=-\lambda\mu_k|g_k|^{p-2}g_k,\quad k\in V\minus\{o\};\label{fr2}\\
\text{boundary  conditions:}&g_o=0.\label{B1}
\end{eqnarray}
where $\{\mu_k: k\in V\}$ is a positive sequence and adopt the convention that $\sum_{i\in \emptyset}f_i=0$ for some sequence $\{f_i\}$ throughout the rest of this paper. If $(\lambda, g)$ with $g\neq0$ is a solution to the eigenvalue problem, then $\lambda$ is called  an p-Laplacian eigenvalue, and $g$ is its eigenfunction. Especially, when $p=2$, the eigenvalue corresponds to the decay rate of birth-death process on trees and $\{\mu_k\}$ is just the invariant measure of birth-death process on trees (see \cite{W-Z}).

Define
$$\aligned
D_p(f)=\sum_{i\in V\setminus\{o\}}\nu_i|f_i-f_{i^*}|^p,\quad f_o=0.
\endaligned$$
Let $(\lambda_p, g)$ be a solution to eigenquation \eqref{fr2} with boundary condition  \eqref{B1}. It is well known that $\lambda_p$ has the following classical variational formula
\begin{equation}\label{f1}\lambda_p=\inf\{D_p(f): \mu(|f|^p)=1, f_o=0\},\end{equation}
We use the ordinary inner product
$$(f, g)=\sum_{k\in V}f_kg_k.$$
Then
$$D_p(g)=(-\Omega_pg, g).$$

Actually, for functions $f$ and $g$ with $f_o=g_o=0$, we have
$$\aligned
(-\Omega_pf, g)&=-\sum_{i\in V}\sum_{j\in J(i)}\nu_j|f_j-f_i|^{p-2} (f_j-f_i)g_i-\sum_{i\in V}\nu_i|f_{i^*}-f_i|^{p-2}(f_{i^*}-f_i)g_i\endaligned$$
By exchanging the order of sums, the formula equals to
$$-\sum_{j\in V\minus\{o\}}\sum_{i=j^*}\nu_j|f_j-f_i|^{p-2} (f_j-f_i)g_i-\sum_{i\in V}\nu_i|f_{i^*}-f_i|^{p-2}(f_{i^*}-f_i)g_i.$$
By $g_o=0$, we have
$$\aligned
(-\Omega_pf, g)&=\sum_{j\in V\minus\{o\}}\nu_j|f_j-f_{j^*}|^{p-2} (f_{j^*}-f_j)(g_{j^*}-g_j).
\endaligned$$
Then the assertion holds by letting $f=g$.

Define
$$\mathscr{D}(D)=\{f: f \text{ is a real function defined on } V,  f_o=0, D_p(f)<\infty\}.$$ Formula \eqref{f1} can be rewritten as the following weighted Hardy inequality:
$$\mu(|f|^p)\leqslant AD_p(f),\quad f\in\mathscr{D}(D),$$
with the optimal constant $A=\lambda_p^{-1}$.
This explains the relationship between the p-Laplacian eigenvalues and the optimal constant of Hardy inequality.

  For a tree $T$,  denote by  $N$ ($N\leqslant\infty$) the maximal level of tree $T$ and $T_i$ ($i$ is included) is a subtree of tree $T$ with $i$  as root. Let  $$\Lambda_i=\{k\in V: |k|=i\},\qquad  i\in\mathbb{Z}^+$$ be the set of elements in the $i^{\text{th}}$ level of the tree. It is clear that $\lambda_p>0$ if $N<\infty$ (otherwise, $\Omega_p g(i)=0$. By letting $i\in \Lambda_N$ in \eqref{fr2}, we have $g_i=g_{i^*}$ for $i\in \Lambda_N$. By the induction, we have $g_i=g_o=0$ for $i\in V$, which is a contraction to $g\neq0$).

  It is easy to see that $\lambda_p=0$ provided $\sum_{k\in V}\mu_k=\infty$ by letting $f_i=1$ for $i\in V\minus\{o\}$ and $f_o=0$ in \eqref{f1}. Therefore, we always assume that $\sum_{k\in V}\mu_k<\infty$. Without loss of generality, we also assume that the root $o$ has a child in the paper.

  We mention that the methods used in this paper are mainly similar to that in \cite{r1}, except one of the key proof of Lemma \ref{Eig} below, in which the monotone of eigenfunction is proved for $p\geqslant2$. Whether Lemma \ref{Eig} still holds for $p\in[1,2)$ or not is still open for us which lead to that some equalities are uncertain in Theorem \ref{th1} below.

The paper is organized as follows. In Section 2, we present the main results, including the monotone of eigenfunction, three kinds of variational formulas for p-Laplacian eigenvalue and its applications (a quantitative estimates  of the p-Laplacian eigenvalue). One example is presented at the end of Section 2. The sketch proofs of the main results are presented in Section 3.

 \section{Main results}
To state our results, we need some  notations.  Let $\mathscr{P}(i)$  be the set of all the vertexes (the root $o$ is excluded) in the unique simple path from $i\in V\minus\{o\}$ to the root and $V_i$ the set of vertexes of subtree $T_i$ for some $i\in V\minus\{o\}$. For $p>1$, let $\hat{p}$ be its conjugate number (i.e., ${1}/{p}+{1}/{\hat{p}}=1$). For $i\in V\setminus\{o\}$, define $\hat{\nu}_j=\nu_j^{1-\hat{p}}$, three operators which are parallel to those introduced in \cite{CWZ1}, as follows:
$$\aligned
&I_i(f)=\frac{1}{\nu_i(f_i-f_{i^*})^{p-1}}\sum_{j\in V_i}\mu_jf_j^{p-1} \qquad \text{(single summation form)},\\
&I\!I_i(f)=\frac{1}{f_i^{p-1}}\bigg[\sum_{k\in \scr{P}(i)}{\hat{\nu}_k}\bigg(\sum_{j\in V_k}\mu_jf_j^{p-1}\bigg)^{\hat{p}-1}\bigg]^{p-1} \qquad \mbox{(double summation form)},\\
&R_i(w)=\mu_i^{-1}\big[\nu_i(1-w_i^{-1})^{p-1}-\sum_{j\in J(i)}\nu_j(w_j -1)^{p-1}\big]\qquad\mbox{(difference form)}.\endaligned$$
Similar operators were  initially  introduced in \cite{Chen1996,Chen2001, r1} respectively for birth-death process in dimension one. We adopt the convention that $1/0=\infty$ and $1/\infty=0$ throughout the paper. To study the lower estimates of p-Laplacian eigenfunction, based on the properties of eigenfunction presented in Lemma \ref{Eig} below, the domains of the three operators are defined respectively as follows:
$$\aligned
\mathscr{F}_I&=\{f: f_o=0, f_i>f_{i^*}\text{ for } i\in V\minus\{o\}\},\endaligned$$
$$\aligned\mathscr{F}_{I\!I}&=\big\{f: f_o=0, f>0 \text{ on } V\minus\{o\} \big\},\endaligned$$
$$\aligned\mathscr{W}&=\big\{w: w>1,  w_o=\infty\}.
\endaligned$$
For the upper bounds,  some modifications are needed to avoid non-summable problem, as shown below.
$$\aligned
\mathscr{{\widetilde F}}_{I}&\!=\!\big\{f\geqslant0: f_o\!=\!0, \exists  k\in V\minus\{o\} \text{ such that } f_i>f_{i^*}\text{ for }i\in\mathscr{P}(k), \text{ and }f_i\!=\!f_{i^*} \text{ for } |i|\!\geqslant\! |k|\big\},\endaligned$$
$$\aligned\mathscr{{\widetilde F}}_{I\!I}&=\{f\geqslant 0: f_o=0, f\neq0, \exists1\leqslant n<N+1 \text{ such that } f_i=f_{i^*} \text{ for } |i|\geqslant n+1\},\endaligned$$
$$\aligned\mathscr{\widetilde{W}}&=\bigcup_{m:\,1\leqslant m<N+1}\bigg\{w: w_o\!=\!\infty, w_i>1\text{ and }\sum_{j\in J(i)}\nu_j(w_j-1)^{p-1}<\nu_i(1-w_i^{-1})^{p-1}\text{ for }
 \\
 &\hskip1.5cm |i|\leqslant m, \text{ and }w_i\!\!=\!1\text{ for } |i|\geqslant m+1\bigg\}.
\endaligned$$
In some extent, these functions are imitated of eigenfunctions of $\lambda_p$. To avoid the trivial estimates,  we need a modified form of  $R$, denote $\widetilde R$, acting on $\widetilde{\scr{W}}$ by replacing $\mu_i$ with $\tilde{\mu}_i=\sum_{j\in V_i}\mu_j$ in $R_i(w)$ when $|i|=m$,
where $m$ is the same one in $\widetilde{\scr{W}}$.  Besides, $\widetilde R$ is also used  when operating the approximating procedure (at this time, $\mu_i$ is replaced with $\tilde \mu_{i}$ for each $i\in V$, see Step 4 in the proof of Theorem \ref{th1} below).
Here and in what follows, the superscript ``$\widetilde{\quad}$'' means modified. The  set below is also needed.
$$\widetilde{\scr{F}}_{I\!I}'=\{f\geqslant0: f_o=0, f\neq0, fI\!I(f)^{\hat{p}-1}\in L^p(\mu)\}.$$

The following lemma  presents us an important property of eigenfunction $g$,  providing the basis for the choices  of those test functions sets of operators $I$, $I\!I$ and $R$. More details see the comments before Lemma \ref{L1} below.
\begin{lem}\label{Eig}Let $T$ be a tree (may have infinite vertexes) with vertexes set $V$ and $p\geqslant2$. If $g\in L^p(\mu)$, $g\neq0$ and $(\lambda_p, g)$ is a solution to \eqref{fr2} with boundary condition $g_o=0$, then $g_i>g_{i^*}$  for each $i\in V\minus\{o\}$.
\end{lem}

In  Theorem \ref{th1}  below, ``$\inf\,\sup$'' are
used for the upper bounds of $\lambda_p$, e.g., each test function $f\in\mathscr{F}_I$
produces a upper bound $\sup_{i\in V\minus\{o\}} I_i(f)^{-1}$, so this part
is called variational formula for upper estimates of $\lambda_p$. Dually, the
``$\sup\,\inf$'' are used for the lower estimates of $\lambda_p$. Among them, the ones expressed by  operator $R$ are easiest
to compute in practice, and the ones expressed by $I\!I$
are hardest to compute but provide better estimates.
Because of ``$\inf\,\sup$'', a localizing procedure is used
for the test function to avoid $I(f)\equiv \infty$ for instance,
which is removed out automatically for the ``$\sup\,\inf$'' part.

\begin{thm}\label{th1} The following variational formulas hold for $\lambda_p$ defined by \eqref{f1}.
\begin{itemize}
\item[$(1)$] Single summation  forms
$$
\aligned \sup_{f\in\mathscr{F}_{I}}
\inf_{i\in V\setminus\{o\}}I_i(f)^{-1}\leqslant\lambda_p\leqslant
\inf_{f\in\mathscr{{\widetilde F}}_{I}}\sup_{i\in V\setminus\{o\}}I_i(f)^{-1},
\endaligned$$
\item[$(2)$] Double summation forms
$$
\aligned
\sup_{f\in S(\mathscr{F})}\inf_{i\in V\setminus\{o\}}I\!I_i(f)^{-1}\leqslant\lambda_p&\leqslant\inf_{f\in S(\widetilde{\mathscr{F}})}\sup_{i\in V\setminus\{o\}}I\!I_i(f)^{-1}
\endaligned$$
with $S(\mathscr{F})=\mathscr{F}_{I\!I}$ or $\mathscr{F}_{I}$ and $S(\widetilde{\mathscr{F}})=\widetilde{\mathscr{F}}_{I\!I}$, or $\widetilde{\mathscr{F}}_{I}$,  or $\widetilde{\scr{F}}_{I\!I}'\cup\widetilde{\scr{F}}_{I\!I}$.
\item[$(3)$] Difference forms
$$
\aligned
 \sup_{w\in\mathscr{W}}\inf_{i\in V\setminus\{o\}}R_i(w)\leqslant\lambda_p\leqslant\inf_{w\in\mathscr{\widetilde{W}}}\sup_{i\in V\setminus\{o\}}\widetilde R_i(w).
\endaligned$$
\end{itemize}
The six equalities in the three terms above hold once $p\geqslant2$.
\end{thm}

We write  $\tilde\mu(A)=\sum_{k\in A}\mu_k$ for some measure $\mu$ and set $A$. Then
$$\mu(V_i)=\sum_{k\in V_i}\mu_k,\quad \hat{\nu}\big(\scr{P}(i)\big)=\sum_{k\in \scr{P}(i)}\hat{\nu}_k, \quad i\in V\minus\{o\}.$$
Define
$$\sigma=\sup_{j\in V\minus\{o\}}\mu\big(T_j\big)\hat{\nu}\big(\scr{P}(i)\big).$$
and $$\#(A)=\text{number of elements in the set } A,$$ for some set $A$.
As applications of Theorem $\ref{th1}$, we have the following theorem.
\begin{thm}\label{Basic}For $p\in(1, \infty)$, we have
$$\sigma^{-1}\geqslant\lambda_p\geqslant\bigg[\Big({\hat{p}}^{p-1}\sup_{i\in V\minus\{o\}}(1+(p-1)C_i)\Big)\, \sigma\bigg]^{-1},$$
where
$$C_i=\#(J(i))+\sum_{s\in J(i)}\sum_{k\in V_s}\big(\#(J(k))-1\big),\quad i\in V.$$
\end{thm}
The theorem effectively presents us the quantitative estimates of the p-Laplacian Dirichlet eigenvalue  on a tree with finite vertexes.
 For the degenerated case of the tree (only one branch), the results reduce to that  on half line in \cite{CWZ1}.
\begin{exam}
Let $T$ be a $r(r\geqslant1)$ order homogeneous tree (i.e., $\#(J(i))=r$, $\forall i\in V\minus\{o\}$) with maximal level $N(\leqslant\infty)$ and root $o$, which has a child, i.e., $\#J(o)=1$. Assume that $t\in(0,1/r)$, $\mu_k=t^{|k|}$ and $\nu_k=ak^{|k|}$ ($a>0$) for $k\in V$. For $p\in(1, \infty)$, denote
$$B_p={\hat{p}}^{p-1}\sup_{i\in V\minus\{o\}}(1+(p-1)C_i).$$
We have
$$\sigma^{-1}\geqslant\lambda_p\geqslant(B_p\sigma)^{-1},$$
where
\begin{eqnarray*}
B_p=\left\{
\begin{array}{ll}{\hat{p}}^{p-1}\big[1+(p-1)\big(r^{N}+r-2\big)\big],&r\geqslant2, \\ p{\hat{p}}^{p-1},&r=1
\end{array}
\right.
\end{eqnarray*}
and
$$\aligned\sigma=\frac{1}{a(1-rt)(1-t^{\hat{p}-1})^{p-1}}\sup_{n\in[1, N+1)}\bigg\{\big(1-(rt)^{N-n+1}\big)(1-t^{n(\hat{p}-1)})^{p-1}\bigg\}.\endaligned$$
If $N=\infty$, then
$$\aligned\sigma=\frac{1}{a(1-rt)(1-t^{\hat{p}-1})^{p-1}}.\endaligned$$
\end{exam}

\section{Proofs of the main results}\label{S2} Without loss of generality, we assume that the root $o$ has only one child (i.e., $\#(J(o))=1$), the  level counting  begins from the child of the  root $o$ (i.e., $|o|=0$ and $|J(o)|=1$),
For convenience, we write   $1$ as the unique child of root $o$ in the proofs of Lemma \ref{Eig}, i.e., $J(o)=\{1\}$ and $\mathscr{P}(1)=\{1\}$.

{\bf Proof of Lemma \ref{Eig}}  We prove the theorem by dividing it into two steps as follows.

(1) we prove that $g_o=0\neq g_1$.

If $g_1=0$, then $\Omega_pg(1)=-\lambda_p\mu_1|g_1|^{p-2}g_1=0$, and
$$\aligned\Omega_pg(1)&=\sum_{j\in J(1)}\nu_j|g_j|^{p-2}g_j.\endaligned$$
Therefore,\begin{equation}\label{eigf1}\sum_{j\in J(1)}\nu_j|g_j|^{p-2}g_j=0.\end{equation}
Moreover, $g_j=0$ for $j\in J(1)$, which will be proved as follows.

Let $A=\{j\in J(1): g_j<0\}$, $B=\{j\in J(1): g_j\geqslant0\}$ and $C_0=\{j\in J(1): g_j=0\}$. Then we prove that $A=B\minus C_0=\emptyset$, which is sufficient to show that $A=\emptyset$ by \eqref{eigf1}. We prove that $A=\emptyset$ by making a contradiction. If $A\neq\emptyset$, then define function $\tilde g$ on $T$ satisfying $\tilde g_o=0$, $\tilde g_1=x>0$, and
\begin{eqnarray*}
\tilde g_{i}=
\left\{
\begin{array}{ll}-g_i,&i\in V_A, \\ g_i,&i\in V_B.
\end{array}
\right.
\end{eqnarray*}
where $V_C:=\cup_{i\in C}V_i$ for some set $C$.
Then
$$\aligned
D_p(\tilde g, \tilde g)&=\sum_{j\in V\minus \{o\}}\nu_j|\tilde g_j-\tilde g_{j^*}|^p\\
&=\nu_1|\tilde g_1-\tilde g_o|^p+\sum_{j\in A}\nu_j|\tilde g_j-\tilde g_{j^*}|^p+\sum_{j\in B}\nu_j|\tilde g_j-\tilde g_{j^*}|^p+\\
&\hskip1cm+\sum_{j\in V_A\minus A}\nu_j|g_j-g_{j^*}|^p+\sum_{j\in V_B\minus B}\nu_j|g_j-g_{j^*}|^p\\
&=\nu_1|x|^p+\sum_{j\in A}\nu_j|-g_j-x|^p+\sum_{j\in B}\nu_j|g_j-x|^p+\\
&\hskip1cm+\sum_{j\in V_A\minus A}\nu_j|g_j-g_{j^*}|^p+\sum_{j\in V_B\minus B}\nu_j|g_j-g_{j^*}|^p\endaligned$$
$$\aligned
D_p(g,g)&=\sum_{j\in V\minus\{o\}}\nu_j|g_j-g_{j^*}|^p\\
&=\sum_{j\in A}\nu_j|g_j|^p+\sum_{j\in B}\nu_j|g_j|^p+\sum_{j\in V_{A\cup B}\minus (A\cup B)}\nu_j|g_j-g_{j^*}|^p\\
&\quad(\text{by } g_1=g_o=0).
\endaligned$$
Therefore,
$$\aligned
&D_p(\tilde g, \tilde g)=D_p(g,g)+\nu_1|x|^p+\sum_{j\in A}\nu_j\big(|g_j+x|^p-|g_j|^p\big)+\sum_{j\in B}\nu_j\big(|g_j-x|^p-|g_j|^p\big)\\
&=D_p(g,g)+\bigg\{\big(\nu_1+\sum_{j\in C_0}\nu_j\big)|x|^p+\sum_{j\in A}\nu_j\big(|g_j+x|^p-|g_j|^p\big)+\sum_{j\in B\minus C_0}\nu_j\big(|g_j-x|^p-|g_j|^p\big)\bigg\}\\
&=:D_p(g,g)+\mu_1G(x).\endaligned$$
We will see that there exists $x>0$ such that $G(x)<0$. Indeed,  let $\ell=\min\{|g_j|/2: j\in A\cup (B\minus C_0)\}$. Then $\ell>0$ by $d_1<\infty$. For $p>0$ and $x\in(0, \ell)$, we have
$|g_j+x|^p<|g_j|^p$ if $j\in A$ and $|g_j-x|^p<|g_j|^p$ if $j\in B\minus C_0$. Since
 $pa^{p-1}(b-a)<b^p-a^p<pb^{p-1}(b-a)$ provided $0<a<b$ and $p\geqslant1$, for $j\in A$ and $x\in(0, \ell)$, we have
$$\aligned|g_j|^p-|g_j+x|^p&\geqslant p|g_j+x|^{p-1}(|g_j|-|g_j+x|)\\
&=p|g_j+x|^{p-1}(-g_j-(-g_j-x))\\
&=p|g_j+x|^{p-1}x\\
&\geqslant\inf_{x\in(0,\ell)}\bigg\{\frac{|g_j+x|^{p-1}}{x^{p-2}}\bigg\}px^{p-1}.\endaligned$$
For $j\in B_0\minus C$ and $x\in(0, \ell)$, we have
$$\aligned|g_j|^p-|g_j-x|^p&\geqslant p|g_j-x|^{p-1}(|g_j|-|g_j-x|)\\
&=p|g_j-x|^{p-1}(g_j-(g_j-x))\\
&=p|g_j-x|^{p-1}x\\
&\geqslant\inf_{x\in(0,\ell)}\bigg\{\frac{|g_j-x|^{p-1}}{x^{p-2}}\bigg\}px^{p-1}.\endaligned$$
When $p\geqslant2$, we have
$$\inf_{x\in(0, \ell)}\bigg\{\frac{|g_j+x|^{p-1}}{x^{p-2}}\bigg\}\geqslant\inf_{x\in(0, |g_j|/2)}\bigg\{\frac{|g_j+x|^{p-1}}{x^{p-2}}\bigg\}\geqslant\frac{|g_j|}{2} \quad\text{ for } j\in A,$$ and
$$\inf_{x\in(0, \ell)}\bigg\{\frac{|g_j-x|^{p-1}}{x^{p-2}}\bigg\}\geqslant\inf_{x\in(0, |g_j|/2)}\bigg\{\frac{|g_j-x|^{p-1}}{x^{p-2}}\bigg\}\geqslant\frac{|g_j|}{2}\quad \text{ for } j\in B\minus C_0.$$
So
$$\aligned &\sum_{j\in A}\nu_j\big(|g_j+x|^p-|g_j|^p\big)+\sum_{j\in B\minus C_0}\nu_j\big(|g_j-x|^p-|g_j|^p\big)\\
&\leqslant -\frac{p}{2}x^{p-1}\bigg\{\sum_{j\in A}\nu_j|g_j|+\sum_{j\in B\minus C_0}\nu_j|g_j|\bigg\}\\
&=:-\frac{p}{2}G_0x^{p-1}<0\quad(p\geqslant2)
.\endaligned$$
Hence,
$$\aligned
G(x)<\big(\nu_1+\sum_{j\in C_0}\nu_j\big)x^p-\frac{p}{2}G_0x^{p-1}.
\endaligned$$
Let $0<x<\min\big\{pG_0\big[2\big(\nu_1+\sum_{j\in C_0}\nu_j\big)\big]^{-1}, |g_j|/2, j\in A\cup(B\minus C_0)\big\}$. Then $G(x)<0$. Moreover,
$$\aligned
D_p(\tilde g, \tilde g)\leqslant D_p(g,g),\quad p\geqslant2.\endaligned$$
Since
$$\mu(|\tilde g|^p)=\sum_{j\in V}\mu_j\tilde g_j=\mu_1x^p+\mu(|g|^p)>\mu(|g|^p),$$
and $g\in L^p(\mu)$, we have
$$\frac{D_p(\tilde g, \tilde g)}{\mu(|\tilde g|^p)}<\frac{D_p(g, g)}{\mu(|g|^p)}\leqslant \lambda_p,$$
which is a contradiction with \eqref{f1}. Therefore, $A=\emptyset$, and $g_j=0$ for $j\in J(1)$. By the induction, we have $g_i=0$ for $i\in V\minus0$. Hence, we must have $g_1\neq0$.

(2) We prove that the eigenfunction  satisfies $g_{i^*}<g_i$ for $i\in V\minus\{o\}$.

We prove the result by making a contradiction. Since $g_1\neq g_o=0$, without loss  of generality, assume that $g_1>0=g_o$(otherwise, replace $g$ by $-g$, which is also an eigenfunction of $\lambda_p$). If there exists $a\in V\minus\{o\}$ satisfying $0=g_o<g_1<\cdots<g_a\geqslant g_b$ for some $b\in J(a)$ ($\mathscr{P}(b)=\{1,\cdots, a,b\}$ and their levels satisfy $|o|\leqslant|1|\leqslant\cdots\leqslant|a|\leqslant|b|$), then set
\begin{eqnarray*}
\overline g_{i}=
\left\{
\begin{array}{ll}g_i,&i\notin V_b, \\ g_a,&i\in V_b.
\end{array}
\right.
\end{eqnarray*}
We have
\begin{eqnarray*}
\Omega_p\overline g(k)&\hskip-2.5cm=\sum_{j\in J(k)}\nu_j|\overline g_j-\overline g_k|^{p-2}(\overline g_j-\overline g_k)+\nu_k|\overline g_{k^*}-\overline g_k|^{p-2}(\overline g_{k^*}-\overline g_k)\\
&\hskip-0.3cm=
\left\{
\begin{array}{lll}0, \qquad\qquad\,\, k\in V_b;\\ \Omega_pg(k), \qquad k\notin V_b,\, k\neq a;\\
\\\sum_{j\in J(a),j\neq b}\nu_j|g_j-g_a|^{p-2}(g_j-g_a)+ \nu_a|g_{a^*}-g_a|^{p-2}(g_{a^*}-g_a), \quad k=a,
\end{array}
\right.\\
&\hskip-4cm=
\left\{
\begin{array}{lll}0, \qquad\qquad k\in V_b;\\ \Omega_pg(k), \qquad k\notin V_b,\, k\neq a;\\
\\ \Omega_pg(a)-\nu_b|g_b-g_a|^{p-2}(g_b-g_a), \qquad k=a,
\end{array}
\right.
\end{eqnarray*}
and
$$\aligned
D_p(\overline g, \overline g)&=(-\Omega_p\overline g, \overline g)_{\mu}=-\sum_{k\in V\minus0}\mu_k\overline g_k\Omega_p\overline g(k)\\
&=-\sum_{k\in V_b}\mu_k\overline g_k\Omega_p\overline g(k)-\sum_{k\notin V_b, k\neq a}\mu_k\overline g_k\Omega_p\overline g(k)-\mu_a\overline g_a\Omega_p\overline g(a)\\
&=-\sum_{k \notin V_b, k\neq a}\mu_k g_k\Omega_pg(k)-\mu_a\overline g_a\Omega_p\overline g(a).
\endaligned$$
By assumption $g_a\geqslant g_b$, we have
 $$\Omega_p\overline g(a)=\Omega_pg(a)-\nu_b|g_b-g_a|^{p-2}(g_b-g_a)\geqslant \Omega_pg(a).$$
 Moreover,
 $$\aligned D_p(\overline g, \overline g)&\leqslant -\sum_{k\notin V_b, k\neq a}\mu_k g_k\Omega_pg(k)-\mu_ag_a \Omega_pg(a)\\
 &=-\sum_{k\notin V_b}\mu_k g_k\Omega_pg(k)\\&=\lambda_p\sum_{k\notin V_b}\mu_k|g_k|^p.\endaligned$$
Since $b\notin\mathscr{P}(1)$, by definition of $\lambda_p$, we have
$$\lambda_p\leqslant\frac{D_p(\overline g, \overline g)}{\mu(|\overline g|^p)}\leqslant\frac{\lambda_p\sum_{k\notin V_b}\mu_k|g_k|^p}{\sum_{k\notin V_b}\mu_k|g_k|^p+\sum_{k\in V_b}\mu_k|g_a|^p}<\lambda_p$$ once $\lambda_p>0$,
which is a contradiction. Therefore, $g_b>g_a$ for each $b\in J(a)$.\qquad$\Box$

Obviously, on the setting of  a finite tree $T$, the eigenfunction $g$ of the p-Laplacian Dirichlet eigenvalue satisfies $g_i>g_{i^*}$ for every $i\in V$.
 Before moving further, we introduce a general equation and discuss the origin of operators used in Theorem \ref{th1}. Recall that $\Lambda_m=\{i: |i|=m\}$ and $N$  is the maximal level of tree $T$. Define $$V(n)=\cup_{m=0}^n \Lambda_m.$$ Consider
$$\text{Poisson equation}: \qquad\Omega_p g(i)=-\mu_i|f_i|^{p-2}f_i,\qquad i\in V\minus\{o\}.$$
By multiplying $\mu_i$ on  both sides of the equation and making summation with respect to $i\in V_k\cap V(n)$ for some $k\in V\minus\{o\}$ with $|k|\leqslant n$, it is easy to check that
\begin{equation}\label{f3-I}\sum_{j\in \Lambda_{n+1}\cap V_k}\nu_j|g_{j^*}-g_j|^{p-2}(g_{j^*}-g_j)-\nu_k|g_k-g_{k^*}|^{p-2}(g_{k^*}-g_k)=\sum_{j\in V_k\cap V(n)}\mu_j |f_j|^{p-2}f_j,\quad |k|\leqslant n.\end{equation}
 If $\lim_{n\to N}\sum_{j\in \Lambda_{n+1}\cap V_k}\nu_j|g_{j^*}-g_j|^{p-2}(g_{j^*}-g_j)=0$ (which is obvious for $N<\infty$),
 then we obtain the form of the operator $I$ by letting $n\to N$ and $f=\lambda_p^{\hat{p}-1}g$ in \eqref{f3-I}.  Moreover, if $g_o=0$ and $g_i>g_{i^*}$ for $i\in V\setminus\{o\}$, then
$$g_i=\sum_{k\in\scr{P}(i)}\hat{\nu}_k\bigg(\sum_{j\in V_k}\mu_j|f_j|^{p-2}f_j\bigg)^{\hat{p}-1}.$$
This explains where the operator $I\!I$ comes from. Similarly, from the eigenequation \eqref{fr2}, we obtain the operator $R$ by letting $w_i=g_i/g_{i^*}$. The eigenequation is a ``bridge'' among these operators.
 Let
$$\tilde\lambda_p=\inf\{D_p(f): \mu(|f|^p)=1, \exists1\leqslant n<N+1 \text{ such that } f_i=f_{i^*} \text{ for } |i|\geqslant n+1\}.$$ If $\sum_{k\in V}\mu_k<\infty$, then $\lambda_p=\tilde\lambda_p$ as will be seen  in Lemma \ref{L1} below.
To this end,  define
$$\lambda_p^{(m)}=\inf\big\{D_p(f): \mu(|f|^p)=1, f_i=f_{i^*} \text{ for }\; |i|\geqslant m+1\big\},\qquad 1\leqslant m<N+1.$$
Let
$$\tilde\mu_i=\mu_i,\quad \tilde\nu_i=\nu_i\qquad\text{ for }\quad |i|\leqslant m-1\text{ and } |j| \leqslant m-1;$$
$$\tilde\mu_i=\sum_{j\in V_i}\mu_j,\quad
\tilde\nu_i=\nu_i\qquad\text{ for }\quad |i|=m.$$ For $f$ with $f_i=f_{i^*}$ for $|i|\geqslant m+1$, we have
$$D_p(f)=\sum_{i\in V(m)\minus\{o\}}\tilde\nu_i|f_i-f_{i^*}|^p=:\widetilde D_p(f),\qquad \mu(|f|^p)=\sum_{i\in V(m)}\tilde\mu_i |f_i|^p=:\tilde\mu(|f|^p).$$
 So  $\lambda_p^{(m)}$ is p-Laplacian eigenvalue of the local Dirichlet form $\big(\widetilde D, \scr{D}(\widetilde D)\big)$ with state space $T(m)$, which is a finite tree with maximal level $m$ and coincides with tree $T$ restricted to the first $m-1$ levels.

 This following lemma presents us an approximating procedure, which guarantees that some properties hold obviously once that hold for finite cases(see  Step 4 in proof of Theorem \ref{th1} below).  For simplicity, we use ``iff'' to denote ``if and only if''  and $\uparrow$(resp. $\downarrow$) to denote increasing and decreasing
  throughout the paper.
\begin{lem}\label{L1} Assume that $\sum_{k\in V}\mu_k<\infty$(i.e., $\mu(T)<\infty$).
We have $\lambda_p=\tilde\lambda_p$ and $\lambda_p^{(n)}\downarrow \lambda_p$ as $n\to N$.
 \end{lem}
 {\noindent{\it Proof}}\quad By definition of $\lambda_p$, for any $\varepsilon>0$, there exists $f$ such that
${D_p(f)}\big/{\mu(|f|^p)}\leqslant \lambda_p+\varepsilon$.
Construct $f^{(n)}$ such that $f_i^{(n)}=f_i$ for $|i|\leqslant n$ and $f_i^{(n)}=f_{i^*}$ for $|i|\geqslant n+1$. Since $\sum_{k\in V}\mu_k<\infty$, we have
$$\aligned
&D_p(f^{(n)})=\sum_{i\in V\minus\{o\}}\nu_i|f_i^{(n)}-f_{i^*}^{(n)}|^p=\sum_{i\in V(n)\minus\{o\}}\nu_i|f_i-f_{i^*}|^p\uparrow D_p(f),\quad n\to N,\\
&\mu({|f^{(n)}|}^p)=\sum_{i\in V(n)\minus\{o\}}\mu_i|f_i|^p+\sum_{i\in \Lambda_{n+1}}\mu(V_i)|f_{i^*}|^p\uparrow \mu(|f|^p)\quad n\to N.\endaligned$$
By definitions of $\lambda_p$, $\tilde\lambda_p$ and $\lambda_p^{(n)}$, the required assertion holds. $\qquad\Box$

 Using the similar methods introduced in \cite{r1}, there are not much difficulties to complete the proof of Theorem \ref{th1}. Therefore, we will present more details of the proofs of Theorem \ref{Basic} but only  some keys for that of  Theorem \ref{th1} in the following.

{\noindent\it\bf Proof of Theorem $\ref{th1}$}\quad
We adopt the following circle to prove the upper bounds of $\lambda_p$.
$$
\aligned
\lambda_p&\leqslant \inf_{f\in\widetilde{\scr{F}}_{I\!I}'\cup\widetilde{\scr{F}}_{I\!I}}\sup_{i\in V\minus\{o\}}I\!I_i(f)^{-1}\\
&\leqslant \inf_{f\in\widetilde{\scr{F}}_{I\!I}}\sup_{i\in V\minus\{o\}}I\!I_i(f)^{-1}=\inf_{f\in\widetilde{\scr{F}}_{I}
}\sup_{i\in V\minus\{o\}}I\!I_i(f)^{-1}=\inf_{f\in\widetilde{\scr{F}}_{I}
}\sup_{i\in V\minus\{o\}}I_i(f)^{-1}\\
&\leqslant \inf_{w\in\widetilde{\scr{W}}}\sup_{i\in V\minus\{o\}}
\widetilde{R_i}(w)\leqslant \lambda_p.
\endaligned$$
The second inequality above is clear and  and the remainders are proved by several steps as follows.

Step 1\quad Prove that $\lambda_p\leqslant \inf_{f\in\widetilde{\scr{F}}_{I\!I}'\cup\widetilde{\scr{F}}_{I\!I}}\sup_{i\in V\minus\{o\}}I\!I_i(f)^{-1}$.

For $f\in \widetilde{\scr{F}}_{I\!I}$, there exists $n\in E$ such that $f_i=f_{i^*}$ for $|i|\geqslant n+1$. Let $$g_i=\sum_{k\in\mathscr{P}(i)}\bigg(\frac{1}{\nu_k}\sum_{j\in V_k}\mu_jf_j^{p-1}\bigg)^{\hat{p}-1}, \quad |i|\leqslant n$$ and $g_i=g_{i^*}$ for $|i|\geqslant n+1$. Then $g\in L^p(\mu)$, $g_i=f_iI\!I_i(f)^{\hat{p}-1}$ for $|i|\leqslant n$ and $f_i\neq0$. Moreover,
$$g_i-g_{i^*}=\bigg(\frac{1}{\nu_i}\sum_{j\in V_i}\mu_jf_j^{p-1} \mathbf{1}_{\{i: |i|\leqslant n\}}\bigg)^{\hat{p}-1}.$$
Inserting this term into $D_p(g)$, we have
$$\aligned D_p(g)&=\sum_{j\in V\minus\{o\}}(g_j-g_{j^*})
\sum_{k\in V_j}\mu_kf_k^{p-1} \mathbf{1}_{\{j: |j|\leqslant n\}}\\
&=\sum_{k\in V\minus\{o\}}\mu_kf_k^{p-1}\sum_{j\in\scr{P}(k)}\mathbf{1}_{\{j: |j|\leqslant n\}}(g_j-g_{j^*})\quad(\text{since } k\in V_j \text{ iff } j\in\scr{P}(k))\\
&=\sum_{k\in V\minus\{o\}}\mu_kf_k^{p-1}g_k\quad(\text{ since }g_i=g_{i^*} \text{ for } |i|\geqslant n+1).\endaligned$$
Since $g\in L^p(\mu)$, we further obtain
$$D_p(g)\leqslant \sum_{k\in V\minus\{o\}}\mu_k |g_k|^p \sup_{k\in  V\minus\{o\}} \bigg(\frac{f_k}{g_k}\bigg)^{p-1}\leqslant \mu(|g|^p) \sup_{k\in  V\minus\{o\}}I\!I_k(f)^{-1}.$$
Hence,
$$\lambda_p\leqslant\frac{D_p(g)}{\mu(|g|^p)}\leqslant\sup_{k\in  V\minus\{o\}}I\!I_k(f)^{-1}.$$
The inequality also holds for $f\in\widetilde{\scr{F}}_{I\!I}'$ since the key point in its proof is $g=fI\!I(f)\in L^p(\mu)$, which also holds for $f\in\widetilde{\scr{F}}_{I\!I}'$. So the required assertion holds.

 Step 2\quad Prove that $$\inf_{f\in\widetilde{\scr{F}}_{I\!I}}\sup_{i\in V\minus\{o\}}I\!I_i(f)^{-1}=\inf_{f\in\widetilde{\scr{F}}_{I}
}\sup_{i\in V\minus\{o\}}I\!I_i(f)^{-1}=\inf_{f\in\widetilde{\scr{F}}_{I}
}\sup_{i\in V\minus\{o\}}I_i(f)^{-1}.$$

(a)\quad We first prove that
$$\inf_{f\in\widetilde{\scr{F}}_{I\!I}}\sup_{i\in V\minus\{o\}}I\!I_i(f)^{-1}\leqslant\inf_{f\in\widetilde{\scr{F}}_{I}
}\sup_{i\in V\minus\{o\}}I\!I_i(f)^{-1}\leqslant\inf_{f\in\widetilde{\scr{F}}_{I}
}\sup_{i\in V\minus\{o\}}I_i(f)^{-1}.$$
Since $\widetilde{\scr{F}}_{I}
\subset\widetilde{\scr{F}}_{I\!I}
$, the first inequality is clear. For $f\in\widetilde{\scr{F}}_{I}
$, there exists $1\leqslant n<N+1$ such that $f_i=f_{i^*}$ for $|i|\geqslant n+1$ and $f_i>f_{i^*}$ for $|i|\leqslant n$.
Since $f_i=\sum_{k\in\scr{P}(i)}(f_k-f_{k^*})$ for $|i|\leqslant n$, inserting this term into the denominator of $I\!I_i(f)$ and using the proportional property, we have
$$\aligned\inf_{i\in V\minus\{o\}}I\!I_i(f)&=\inf_{i\in V(n)\minus\{o\}}I\!I_i(f)\geqslant\inf_{i\in V\minus\{o\}}I_i(f). \endaligned$$
and the required assertion holds since $f\in\widetilde{\scr{F}}_{I}$ is arbitrary.

(b)\quad Prove the equality.

For $f\in\widetilde{\scr{F}}_{I\!I}$, $\exists\, n\in[1, N+1)$ such that $f_i=f_{i^*}$ for $|i|\geqslant n+1$ and $f\neq0$. Let $$g_i=\sum_{k\in\mathscr{P}(i)}\bigg(\frac{1}{\nu_k}\sum_{j\in V_k}\mu_jf_j^{p-1}\bigg)^{\hat{p}-1}, \quad 0<|i|\leqslant n,$$ $g_o=0$ and $g_i=g_{i^*}$ for $|i|\geqslant n+1$. Then $g\in\widetilde{\scr{F}}_{I}$ and
$$g_i-g_{i^*}=\bigg(\frac{1}{\nu_i}\sum_{j\in V_i}\mu_j f_j^{p-1}\bigg)^{\hat{p}-1},\qquad 0<|i|\leqslant n.$$
Moreover,
$$\nu_i(g_i-g_{i^*})^{p-1}\leqslant \sum_{j\in V_i}\mu_j g_j^{p-1}\sup_{j\in V_i}\bigg(\frac{f_j}{g_j}\bigg)^{p-1}=\sum_{j\in V_i}\mu_j g_j^{p-1}\sup_{j\in V_i}I\!I_i(f)^{-1},\qquad i\in V\minus\{o\}.$$
Hence,
$$\sup_{k\in V\minus\{o\}}I_k(g)^{-1}\leqslant\sup_{k\in V\minus\{o\}}I\!I_k(f)^{-1}.$$
Then the assertion follows by making the infimum over $\widetilde{\scr{F}}_{I}$ first and then the infimum over $\widetilde{\scr{F}}_{I\!I}$.

Step 3\quad Prove that $\inf_{f\in\widetilde{\scr{F}}_{I\!I}}\sup_{i\in V\minus\{o\}}I\!I_i(f)^{-1}\leqslant\inf_{w\in\widetilde{\scr{W}}}\sup_{i\in V\minus\{o\}}
\widetilde{R_i}(w)$.

First, we change the form of $\widetilde R$. For $w\in \widetilde{\scr{W}}$ with $w_i=1$ for $|i|\geqslant m+1$, let $g$ be a positive function on $V\minus\{o\}$ with $g_o=0$  such that $w_i=g_i/g_{i^*}$. Applying  Lemma \ref{Eig} to the finite tree $T(m)$, we have  $g_i>g_{i^*}$ for $|i|\leqslant m$ and $g_i=g_{i^*}$ for $|i|\geqslant m+1$. Since
$$\sum_{j\in J(i)}\nu_j(w_j-1)^{p-1}<\nu_i(1-w_i^{-1})^{p-1}\quad\text{ for }\quad |i|\leqslant m,$$
 we have
$\widetilde R_i(w)={-\widetilde\Omega_p g(i)}\big/\mu_i{g_i}^{p-1}>0$ for $|i|\leqslant m$ and $\widetilde R_i(w)=0$ for $|i|\geqslant m+1$,
where $\widetilde\Omega$ is a change form of $\Omega$ with $\mu_i$ replaced by $\tilde \mu_i$ for $|i|\leqslant m$.

Now, we come back to the main assertion.  For $w\in \widetilde{\scr{W}}$ with $w_i=1$ for $|i|\geqslant m+1$, let $g$ be the function mentioned above and
\begin{eqnarray*}
f_{i}^{p-1}=
\left\{
\begin{array}{lll}
-\mu_i^{-1}\bigg[\sum_{j\in J(i)}\nu_j(g_j-g_i)^{p-1}+\nu_i(g_i-g_{i^*})^{p-1}\bigg],&|i|\leqslant m-1, \\
\tilde \mu_i^{-1}\nu_i(g_i-g_{i^*})^{p-1},&|i|=m,\\
f_{i^*}^{p-1},&|i|\geqslant m+1.
\end{array}
\right.
\end{eqnarray*}
Then $\mu_if_i^{p-1}=-\widetilde\Omega_p g(i)>0$ for $|i|\leqslant m$. By \eqref{f3-I}, we have
\begin{equation}\label{f4}\nu_k(g_k-g_{k^*})^{p-1}=\sum_{j\in V_k\cap V(m-1)}\mu_j f_j^{p-1}+\sum_{j\in \Lambda_{m}\cap V_k}\nu_j(g_j- g_{j^*})^{p-1},\quad |k|\leqslant m-1.
\end{equation}
 Since
$$\nu_j(g_j-g_{j^*})^{p-1}=\sum_{i\in V_j}\mu_if_j^{p-1}=\sum_{i\in V_j}\mu_if_i^{p-1},\qquad|j|=m,$$
we have
$$\aligned\sum_{j\in \Lambda_{m}\cap V_k}\nu_j(g_j-g_{j^*})^{p-1}&=\sum_{j\in \Lambda_m\cap V_k}\sum_{i\in V_j}\mu_if_i^{p-1}=\sum_{j\in \big(T\minus V(m-1)\big)\cap V_k}\mu_jf_j^{p-1}\\
&=\sum_{j\in V_k}\mu_jf_j^{p-1}-\sum_{j\in V(m-1)\cap V_k}\mu_jf_j^{p-1},\quad |k|\leqslant m.\endaligned$$
Inserting this term into \eqref{f4}, we arrived at
$$\aligned\nu_k(g_k-g_{k^*})^{p-1}&=\sum_{j\in V_k\cap V(m-1)}\mu_j f_j^{p-1}+\sum_{j\in \Lambda_{m}\cap V_k}\nu_j(g_j-g_{j^*})^{p-1}\\
&=\sum_{j\in V_k}\mu_jf_j^{p-1},\quad 0<|k|\leqslant m-1.\endaligned$$
Hence,
$$\nu_k(g_k-g_{k^*})^{p-1}=\sum_{j\in V_k}\mu_jf_j^{p-1}\quad 0<|k|\leqslant m.$$
Moreover,
$$g_i=\sum_{k\in\scr{P}(i)}\bigg(\frac{1}{\nu_k}\sum_{j\in V_k}\mu_jf_j^{p-1}\bigg)^{\hat{p}-1}\,\qquad 0<|i|\leqslant m,$$
and $\widetilde R_i(w)=(f_i/g_i)^{p-1}=I\!I_i(f)^{-1}$ for $0<|i|\leqslant m$.
Since $\widetilde R_i(w)=0$ and $f_i=f_{i^*}$ for $|i|\geqslant m+1$, we obtain
$$\sup_{i\in V\minus\{o\}}\widetilde R_i(w)=\sup_{i\in V\minus\{o\}}I\!I_i(f)^{-1}\geqslant\inf_{f\in\widetilde{\scr{F}}_{I\!I}}\sup_{i\in V\minus\{o\}}I\!I_i(f)^{-1},\qquad w\in\widetilde{\scr{W}},$$
and the required assertion holds.

Step 4\quad Prove that $\inf_{w\in\widetilde{\scr{W}}}\sup_{i\in V\minus\{o\}}\widetilde R_i(w)\leqslant\lambda_p$ once $p\geqslant2$.

 Let $g$ with $g_o=0$ be an eigenfunction of local p-Laplacian eigenvalue $\lambda_p^{(m)}$ and extend $g$ to  $V\minus\{o\}$ by setting $g_i=g_{i^*}$ for  $|i|\geqslant m+1$. Put $w_i=g_i/g_{i^*}$ for $i\in V\minus\{o\}$. Then $w\in\widetilde{\scr{W}}$ provided $p\geqslant2$.

Since $m<\infty$, we have $\widetilde R_i(w)=\lambda_p^{(m)}>0$ for $i\in V(m)\minus\{o\}$ and $\widetilde R_i(w)=0$ for $T\minus V(m)$.
 Therefore,
 $$\aligned
 \lambda_p^{(m)}&=\sup_{i\in V\minus\{o\}}\widetilde R_i(w)\\
 &\geqslant \inf_{w\in\widetilde{\scr{W}}: w_i=1\text{ for }|i|\geqslant m+1}\sup_{i\in V(m)\minus\{o\}}\widetilde R_i(w)\\
 &\geqslant\inf_{w\in\widetilde{\scr{W}}: \exists n\geqslant1 \text{ such that } w_i=1\text{ for } | i |\geqslant n+1}\sup_{i\in V\minus\{o\}}\widetilde R_i(w)\\
 &\geqslant \inf_{w\in\widetilde{\scr{W}}}\sup_{i\in V \minus\{0\}}\widetilde R_i(w).
 \endaligned$$
Letting $m\to N$, the assertion then follows by Lemma \ref{L1}.

By now, we have finished the proof of the estimates of upper $\lambda_p$. Dually,
one can prove the lower estimates without too much difficulty. We ignore the
details here.
\qquad$\Box$

Define $T_{i, j}=T_i\cup T_j$, correspondingly, $V_{i, j}=V_i\cup V_j$.  Then $$V_{J(i)}=\big\{k: s\in J(i) \text{ and } k\in V_s\big\}.$$
Similarly, $$J(V_i)=\{k: s\in V_i\text{ and } k\in J(s)\}.$$
It is obvious that $J(V_i)=V_{J(i)}$. Without loss of generality, we adopt convention that $\mu\big(V_k\big)=0$  if $V_k=\phi$.  We use notation $|J(i)|$ indicating the level of $J(i)$, i.e., $|(J(i))|=|i|+1$. For simplicity, we also write $\varphi_i=\hat{\nu}\big(\scr{P}(i)\big)^{p-1}$ in the proof of Theorem $\ref{Basic}$.

{\noindent\it\bf Proof of Theorem $\ref{Basic}$}\quad First, we prove that $\lambda_p^{-1}\leqslant\big({\hat{p}}^{p-1}\sup_{i\in V\minus\{o\}}\big(1+(p-1)C_i\big)\big)\, \sigma$. By calculation, we have
$$\aligned\sum_{j\in V_i}\mu_jf_j^{p-1}&=\sum_{j\in V_i}f_j^{p-1}\bigg[\mu(V_j)-\sum_{k\in J(j)}\mu(V_k)\bigg]\\
&=\sum_{j\in V_i}\mu(V_j)f_j^{p-1}-\sum_{j\in V_i}\sum_{k\in J(j)}\mu(V_k)f_j^{p-1}\\
&=\sum_{j\in V_i}\mu(V_j)f_j^{p-1}-\sum_{k\in V_{J(i)}}\mu(V_k)f_{k^*}^{p-1}\quad(\text{since } J(V_i)=V_{J(i)})\\
&=\mu(V_i)f_i^{p-1}+\sum_{k\in V_{J(i)}}\mu(V_k)\big(f_k^{p-1}-f_{k^*}^{p-1}\big)\quad\big(\text{since }V_i=\{i\}\cup V_{J(i)}\big).\endaligned$$
Put $f_j=\varphi_j^{1/p}$ for $j\in V\minus\{o\}$. Then
$$\aligned\sum_{j\in V_i}\mu_jf_j^{p-1}&=\mu(V_i)\varphi_i^{1/\hat{p}}+\sum_{k\in V_{J(i)}}\mu(V_k)\big(\varphi_k^{1/\hat{p}}-\varphi_{k^*}^{1/\hat{p}}\big)\\
&\leqslant \sigma\bigg[{\varphi_i}^{-1/p}+\sum_{k\in V_{J(i)}}\frac{1}{\varphi_k}\bigg(\varphi_k^{1/\hat{p}}-\varphi_{k^*}^{1/\hat{p}}\bigg)\bigg].
\endaligned$$
Since $\varphi_k\geqslant\varphi_{k^*}$, we obtain
\begin{equation}\label{ad}
\sum_{k\in V_{J(i)}}\frac{1}{\varphi_k}\bigg(\varphi_k^{1/\hat{p}}-\varphi_{k^*}^{1/\hat{p}}\bigg)\leqslant
(p-1)\sum_{k\in V_{J(i)}}\big(\varphi_{k^*}^{-1/p}-\varphi_k^{-1/p}\big).
\end{equation}
Indeed, it suffices to show that
$$\aligned \varphi_k^{1/\hat{p}}-\varphi_{k^*}^{1/\hat{p}}\leqslant
(p-1){\varphi_k}\big(\varphi_{k^*}^{-1/p}-\varphi_k^{-1/p}\big),
\endaligned$$
or equivalently,
$$\aligned p\varphi_k^{1/\hat{p}}-\varphi_{k^*}^{1/\hat{p}}\leqslant
(p-1)\varphi_k\varphi_{k^*}^{-1/p},
\endaligned$$
i.e. $$\aligned p\varphi_k^{1/\hat{p}}\varphi_{k^*}^{1/p}\leqslant\varphi_{k^*}+(p-1)\varphi_k=p\bigg(\frac{1}{p}\big(\varphi_{k^*}^{1/p}\big)^p+\frac{1}{\hat{p}}\big(\varphi_k^{1/\hat{p}}\big)^{\hat{p}}\bigg),
\endaligned$$
which is obvious by Young's inequality.

Noticing that $V_{J(i)}=J(V_i)$ and $k\in J(j)$ if and only if $k^*=j$, we have
$$\aligned
\sum_{k\in V_{J(i)}}\varphi_{k^*}^{-1/p}= \sum_{k\in J(V_i)}\varphi_{k^*}^{-1/p} &=\sum_{j\in V_{i}}\sum_{k\in J(j)}\varphi_j^{-1/p}=\sum_{j\in V_{i}}\#(J(j))\varphi_j^{-1/p},
\endaligned$$
Inserting the term to the inequality \eqref{ad}, we get
$$\aligned
\sum_{k\in V_{J(i)}}\frac{1}{\varphi_k^{1/p}}\big(\varphi_k^{1/\hat{p}}-\varphi_{k^*}^{1/\hat{p}}\big)
&\leqslant
(p-1)\bigg\{\sum_{j\in V_{i}}\#(J(j))\varphi_j^{-1/p}-\sum_{k\in V_{J(i)}}\varphi_k^{-1/p}\bigg\}\\
&=(p-1)\bigg\{\#(J(i))\varphi_i^{-1/p}+\sum_{k\in V_{J(i)}}\bigg(|J(k)|-1\bigg)\varphi_k^{-1/p}\bigg\}\\
&\leqslant(p-1)\bigg[\#(J(i))+\sum_{k\in V_{J(i)}}\bigg(|J(k)|-1\bigg)\bigg]\varphi_i^{-1/p}\quad\big(\text{since }\varphi_k\geqslant\varphi_{k^*}\big).
\endaligned$$
Hence,
$$\aligned\sum_{j\in V_i}\mu_j\varphi_j^{1/p}&\leqslant\bigg[1+(p-1) \bigg(\#(J(i))+\sum_{s\in J(i)}\sum_{k\in V_s}\bigg(|J(k)|-1\bigg)\bigg)\bigg]\sigma\varphi_i^{-1/p}
\\&=\big(1+(p-1)C_i\big)\sigma\varphi_i^{-1/p}.
\endaligned$$
Since
$$\big(\varphi_i^{\hat{p}-1}-\varphi_{i^*}^{\hat{p}-1}\big)^{p-1}=\frac{1}{\nu_i}\quad\text{ and }\quad \varphi_i^{1/[p(p-1)]}\varphi_{i^*}^{1/p}\leqslant\frac{1}{p}\varphi_i^{\hat{p}-1}+\frac{1}{\hat{p}}\varphi_{i^*}^{\hat{p}-1},$$
we obtain
$$\aligned
I_i(\varphi^{1/p})&=\frac{1}{\nu_i\big(\varphi_i^{1/p}-\varphi_{i^*}^{1/p}\big)^{p-1}}\sum_{j\in V_i}\mu_j\varphi_j^{1/\hat{p}}\\
&\leqslant \big(1+(p-1)C_i\big)\sigma\varphi_i^{-1/p}\bigg(\frac{\varphi_i^{\hat{p}-1}-\varphi_{i^*}^{\hat{p}-1}}{\varphi_i^{1/p}-\varphi_{i^*}^{1/p}}\bigg)^{p-1}\\
&=\big(1+(p-1)C_i\big)\sigma\bigg(\frac{\varphi_i^{\hat{p}-1}-\varphi_{i^*}^{\hat{p}-1}}{\varphi_i^{\hat{p}-1}-\varphi_{i^*}^{1/p}\varphi_i^{1/[p(p-1)]}}\bigg)^{p-1}\\
&\leqslant\big(1+(p-1)C_i\big){\hat{p}}^{p-1}\sigma.
\endaligned$$
It is clear that $\varphi^{1/p}\in\scr{F}_I$, we have
$$\lambda_p^{-1}\leqslant\inf_{f\in\scr{F}_I}\sup_{i\in V\minus\{o\}}I_i(f)\leqslant\sup_{i\in V\minus\{o\}}I_i(\varphi^{1/p})\leqslant\Big({\hat{p}}^{p-1}\sup_{i\in V\minus\{o\}}\big(1+(p-1)C_i\big)\Big)\sigma,$$
by Theorem \ref{th1}\,(1).

Now, we prove that $\lambda_p\leqslant\sigma^{-1}$.
For $i_0\in V\minus\{o\}$, let $f$ be a function such that
\begin{equation*} f_i=
\begin{cases}
    \varphi_i^{\hat{p}-1} &  \text{if }\; i\in \scr{P}(i_0), \\
    \varphi_{i_0}^{\hat{p}-1} &  \text{if}\; i\in V_{i_0},\\
    0   &\;  \text{Others}.
 \end{cases}
 \end{equation*}
Then
$$\aligned
\mu(|f|^p)&=\sum_{i\in V\minus\{o\}}\mu_i|f_i|^p=\sum_{i\in\scr{P}(i_0)}\mu_i\varphi_i^{\hat{p}}+\mu(T_{i_0})\varphi_{i_0}^{\hat{p}}.\endaligned$$
Since $f_i-f_{i^*}=\big(\nu_i\big)^{\hat{p}-1}$ for $i\in\scr{P}(i_0)$ and $f_i-f_{i^*}=0$ for $i\in V\minus\scr{P}(i_0)$.
$$\aligned D_p(f)&=\sum_{i\in V\minus\{o\}}\nu_i|f_i-f_{i^*}|^p=\varphi_{i_0}^{\hat{p}-1}.
\endaligned$$
By \eqref{f1}, we have
$$\lambda_p^{-1}\geqslant\frac{\mu{|f|^p}}{D_p(f)}\geqslant \mu(T_{i_0})\varphi_{i_0},\qquad i_0\in V\minus\{o\}.$$
Then the assertion follows by taking supremum over $V\minus\{o\}$.$\qquad\Box$

\bigskip


\end{document}